\documentclass{article}
\usepackage{amsmath,amssymb,amsfonts,theorem}
\usepackage{graphicx,epsfig,color} %psfrag,xspace
\graphicspath{{./Figures/}}
\newcommand{\until}[1]{\{1,\dots, #1\}}
\newcommand{\subscr}[2]{#1_{\textup{#2}}}

\newcommand{\map}[3]{#1: #2 \rightarrow #3}
\newcommand{\nn}{\nonumber}
\newtheorem{theorem}{Theorem}[section]
\newtheorem{proposition}[theorem]{Proposition}

{\theorembodyfont{\rmfamily} % These ones have roman font in their body
\newtheorem{remark}[theorem]{Remark}

}

\newcommand{\Exp}{\mathbb{E}}
\newcommand{\R}{\mathbb{R}}

\newcommand{\nonnegativeinteger}{\mathbb{Z}_{\ge0}}
\newcommand{\G}{\mathcal{G}}
\newcommand{\beqn}{\begin{eqnarray}}
\newcommand{\eeqn}{\end{eqnarray}}

\newcommand{\argmin}{\mathop{\operatorname{argmin}}}
\newcommand{\1}{\mathbf{1}} % vector of ones
\newcommand{\neigh}{\mathcal{N}}
\newcommand{\reals}{\R}
\newcommand{\realpositive}{\ensuremath{\reals_{>0}}}

\newcommand{\xave}{\subscr{x}{ave}}
\newcommand{\lopt}{\subscr{l}{opt}}
\newcommand{\centerG}{m}
\title{Detection of Gaussian signals via hexagonal sensor networks}
\author{Paolo~Frasca \and Paolo~Mason \and Benedetto~Piccoli %
\thanks{This work has been partially supported by the FIRB~2005 CASHMA research project. %
Paolo~Frasca is with the D.I.I.M.A., Universit\`{a} di Salerno, and with the Istituto per le Applicazioni del Calcolo-C.N.R., Roma, Italy. %
Paolo~Mason is with the Ecole sup\'erieure d'\'electricit\'e, Paris, France. %
Benedetto~Piccoli is with the Istituto per le Applicazioni del Calcolo-C.N.R., Roma, Italy.
}
}

\begin{document}
\maketitle
\thispagestyle{empty}
\pagestyle{empty}
\begin{abstract}
This paper considers a special case of the problem of identifying a static scalar signal, depending on the location, using a planar network of sensors in a distributed fashion. Motivated by the application to monitoring wild-fires spreading and pollutants dispersion, we assume the signal to be Gaussian in space. Using a network of sensors positioned to form a regular hexagonal tessellation, we prove that each node can estimate the parameters of the Gaussian from local measurements. Moreover, we study the sensitivity of these estimates to additive errors affecting the measurements. Finally, we show how a consensus algorithm can be designed to fuse the local estimates into a shared global estimate, effectively compensating the measurement errors.
\end{abstract}

\section{Introduction}
In this paper we consider a distributed parametric estimation problem: we want to identify a static scalar field, depending on the location, using a planar network of sensors in a distributed fashion. We assume the field to be Gaussian in space: the goal of the sensor network is to identify the four parameters describing the Gaussian function.
This assumption can be seen as a simplified model of pollutants dispersion or of wild-fires spreading, which the network is required to monitor.  The objective of this paper is to design a network able to achieve such identification goal, using noisy punctual measurements and local communication. We are going to design the spacial deployment of the sensors, their communication topology, and the algorithms both to infer the parameters from the measurements, and to aggregate local information into a common estimate.

Sensor networks and distributed estimation are a vast and rapidly growing field: among the available literature, we refer the reader to \cite{AS-CF-KHJ-ASV:08,GCC-FA:09,FG-SZ:09,SHD-FB:07n} for some recent research works, and to \cite{MR-RN:04,AG-PRK:06,AG-SKD:08} for some broader discussions.

Among such wide literature, three recent works \cite{KML-IRS-PY-RAF:08,SM:07a,FZ-NEL:09} are especially relevant to us, as they consider the problem of estimating a scalar planar field through point-wise noisy measurements, taken by coordinated moving sensors. Let us discuss their contributions. First, the authors of \cite{KML-IRS-PY-RAF:08} consider the problem of finding the leading coefficients of a suitable series expansion of the field. The proposed solution relies on a distributed version of a filter, based on average consensus.
Second, the forthcoming paper \cite{SM:07a} develops a non parametric interpolation algorithm: the field is defined on a compact subset of the plane, which is apportioned into regions assigned to each agent. Then, each node maintains a representation of its current region. Third, the papers \cite{PO-EF-NEL:04,FZ-NEL:09} contain an articulated framework of cooperative exploration, in which a group of sensors travels in the plane, harvesting information about the field.

With respect to these works, our paper differs on many respects. First, as we are interested in monitoring rather than in exploration, we choose to consider a static network, rather than a group of moving sensors.
More important, we undertake one particular parametric estimation problem, which has not received any specific attention in the literature. This allows to obtain simple and neat results, and to focus on the core problem of fusing the local estimates obtained by the single nodes.

%Moreover, the hexagonal tessellation has been widely studied in the literature, for many applications, see \cite{ESB-GS:03,TA-NEV-MRS-PVK:07}.

\subsection{Problem statement}\label{sec:ProblemStatement}
Let $\map{F}{\R^2}{\R}$ be a Gaussian function
\begin{equation}\label{eq:Gaussian}
F(x_1,x_2)=C_1 e^{-\frac{(x_1 - \centerG_1)^2 +(x_2 - \centerG_2)^2}{C_2}}
\end{equation}
where $C_1$ and $C_2$ are positive constants and $\centerG=(\centerG_1,\centerG_2)$ is a point of the plane. Let us discuss the meaning of these parameters, considering the application of observing a diffusion process (fire, pollution).
\begin{enumerate}
\item The point $\centerG$ is the center of the distribution, that is the center of the fire or the location of the pollutant leak.
\item The positive scalar $C_1$ gives the amplitude of $F$, that is $C_1=F(\centerG_1,\centerG_2)$.
\item The positive scalar $C_2$ gives information about the width of the event. If the distribution is the outcome of a diffusion process, $C_2$ is proportional to the square of the time elapsed from the beginning of the event.
\end{enumerate}
Let there be a network of sensors, referred to as nodes. Each sensor is able to measure the value of $F$ in its own location, and to communicate with its neighbors on the graph with no errors or delays.
We want the network to be able to collectively estimate of the four parameters describing the Gaussian function.
To this goal, we devise a two-phases algorithm.
\begin{enumerate}
\item[1)] The sensors, communicating with their neighbors, construct local estimates of the parameters.
\item[2)] The local estimates are fused into common estimates, using an iterative consensus algorithm.
\end{enumerate}

\subsection{Statement of contributions and outline}
We design the two phases of distributed sensing and fusion separately. In Section~\ref{sec:Network} we discuss the geometric design of the sensor network, and we choose to work with an hexagonal tessellation. In Section~\ref{sec:Estimation} we consider how a single node of such network can build an estimate of the parameters of the Gaussian field, using local measurements, that is measurements taken by itself and its neighbors. Moreover, we study the the sensitivity to measurement errors of the estimates of the parameters, and the way the induced estimation error depends on the distance among sensors.
Afterwards, in Section~\ref{sec:Fusion}, we consider the whole network, and we discuss how local estimates can be fused into a common estimate, through an iterative consensus algorithm. We show that fusion can effectively compensate errors in the measurements. Finally, Section~\ref{sec:outro} contains some concluding remarks.

%%%%%%%%%%%%%%%%%%%%%%%%%%%%%%%%%%%%%%%%%%%%%%%%%%%%%%%%%%%%%%%%%%%%%%%%%%%%%%%%%%%%
\section{Network design: hexagonal tessellation}\label{sec:Network}
In this section we discuss the properties for a network of sensors which are desirable to solve the problem stated in Section~\ref{sec:ProblemStatement}. We choose to work with networks whose nodes and links form (a portion of) a regular {\em hexagonal tessellation} of the plane. Indeed, an hexagonal tessellation enjoys many properties which make it the right candidate for our network.
\begin{enumerate}
\item {\em Equispacing.} In an hexagonal tessellation, all neighboring nodes are $l$ far from each other, if we denote by $l$ the length of the edge of each regular hexagon. If the communication among nodes is wireless, which is a very natural choice for a wide area network to be deployed in an environment, a common communication model is the disk graph: each node can communicate with nodes closer than a threshold $R$. Since in a regular hexagonal tessellation all nodes are $l$ far apart from their neighbors, choosing $l<R$ we satisfy the connectivity constraint.
\item {\em Modularity.} Neglecting border effects, an hexagonal tessellation is vertex-transitive. This implies that the algorithmic design can easily be done off-line, and then be applied to any node.
\item {\em Coverage properties.} Regular tessellations are optimal from the point of view of covering, since they are the Delaunay graph of a centroidal Voronoi partition \cite{FB-JC-SM:09}.
\end{enumerate}
In a distributed framework, every single node is supposed to obtain an estimate of the four parameters $C_1,C_2,\centerG_1,\centerG_2$. We are going to prove in the sequel that, in order to do this, it needs at least four measurements, which can be its own measurement and three measurements from its neighbors. The graph of the hexagonal tessellation is {\em 3-regular} (every node has three neighbors), and then it allows the required information exchange.

The above properties are satisfied not only by hexagonal tessellations, but also by the other regular tessellations of the plane, using triangles or squares. However, hexagons are preferable. Indeed, they require the minimum number of connections, three (instead of six and four, respectively): this can be a desirable property in terms of the communication cost. Moreover, hexagonal tessellations cover a larger area with the same number of nodes, and this is a useful feature, considering that the cost of putting in place a wireless network is proportional to the number of nodes. Indeed elementary formulas imply that, being $N$ the number of nodes, the covered area is $\frac{3\sqrt{3}}{4} l^2 N$ for the hexagonal tessellation, while it is $\frac{\sqrt{3}}{2} l^2 N,$ and $l^2 N,$ for the triangular and square tessellations, respectively.

\begin{remark}
The regularity properties of the network are lost at the borders. For this reason, estimation is going to be performed at inner nodes only: such limitation is intuitively going to be less important for large networks.
\end{remark}
%, since the number of external nodes is proportional to $\sqrt{N}.$

%%%%%%%%%%%%%%%%%%%%%%%%%%%%%%%%%%%%%%%%%%%%%%%%%%%%%%%%%%%%%%%%%%%%%%%%%%%%%%%%%%%%
\section{Estimation}\label{sec:Estimation}
\begin{figure}\begin{center}
\includegraphics[width=.6\columnwidth]{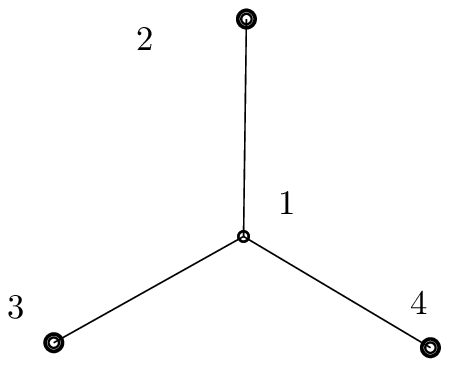}
\caption{Local view of the neighborhood of a generic node, with the local labeling used in Proposition~\ref{prop:inverting}.}
\label{fig:local}
\end{center}\end{figure}

In this section, we show how an inner node of the hexagonal tessellation graph can estimate the four parameters of the Gausssian, using local information.
Since all edges are of length $l$, up to isometries we can assume the four sensors to be located in
$$
\left(
\begin{array}{c}
0 \\ 0
\end{array}
\right),\
\left(
\begin{array}{c}
0 \\ l
\end{array}
\right),\
\left(
\begin{array}{c}
-\frac{\sqrt{3}}2 l \\ -\frac{l}2
\end{array}
\right),\
\left(
\begin{array}{c}
\frac{\sqrt{3}}2 l \\ -\frac{l}2
\end{array}
\right),
$$
and labeled as in Figure~\ref{fig:local}.

\begin{proposition}\label{prop:inverting}
Let $\mu_1,\ \mu_2,\ \mu_3,\ \mu_4$ be the measurements taken by the four sensors, labeled as in Figure~\ref{fig:local}. Then the four parameters in \eqref{eq:Gaussian} are given by the following formulas:
\begin{align*}
C_2=&\frac{3 l^2}{\log (\frac{\mu_1^3}{\mu_2\mu_3\mu_4})};\\%\label{eq:C2}\\
\centerG_1=&\frac{C_2}{2 l \sqrt{3}} \log{\frac{\mu_4}{\mu_3}};\\
\centerG_2=&\frac{C_2}{6 l} \log{\frac{\mu_2^2}{\mu_3\mu_4}};\\
C_1=& \mu_1 e^{\frac{\centerG_1^2+\centerG_2^2}{C_2}}.
\end{align*}
\end{proposition}
\begin{proof}
Let $\map{\Phi}{\realpositive\times\realpositive\times\reals\times\reals}{\realpositive^4}$ express the four measurements as functions of the parameters.
\begin{align}\label{eq:def-Phi}\Phi(&C_1,C_2,\centerG_1,\centerG_2)
=&\Big(F(0,0),F(0,l),F\Big(-\frac{\sqrt{3}}2 l, -\frac{l}2\Big),F\Big(\frac{\sqrt{3}}2 l,-\frac{l}2\Big)\Big)\end{align}
Then the problem consists in inverting $\Phi$, that is solving the algebraic system
$$
\left\{\begin{array}{lcl}
C_1 e^{-\frac{\centerG_1^2 +\centerG_2^2}{C_2}}&=&\mu_1 \\
C_1 e^{-\frac{\centerG_1^2 +(l-\centerG_2)^2}{C_2}}&=&\mu_2 \\
C_1 e^{-\frac{(\frac{\sqrt{3}}2 l+\centerG_1)^2 +(\frac{l}2+\centerG_2)^2}{C_2}}&=&\mu_3 \\
C_1 e^{-\frac{(\frac{\sqrt{3}}2 l-\centerG_1)^2 +(\frac{l}2+\centerG_2)^2}{C_2}}&=&\mu_4.
\end{array}\right.
$$
This implies $\frac{\mu_1^3}{\mu_2\mu_3\mu_4}=e^{\frac{3 l^2}{C_2}},$ which requires that $\mu_2\mu_3\mu_4<\mu_1^3$.
%$$(\mu_1,\mu_2,\mu_3,\mu_4)\in \{(\mu_1,\mu_2,\mu_3,\mu_4)\in\realpositive^4\ :\ \mu_2\mu_3\mu_4<\mu_1^3\}.$$
If this holds, we can deduce $C_2=\frac{3 l^2}{\log (\frac{\mu_1^3}{\mu_2\mu_3\mu_4})}.$
Moreover, $\frac{\mu_4}{\mu_3}=e^{\frac{2 l \sqrt{3}\centerG_1}{C_2}}\,,\quad \frac{\mu_2^2}{\mu_3\mu_4}=e^{\frac{4l\centerG_2}{C_2}},$
so that
$$\centerG_1=\frac{C_2}{2 l \sqrt{3}} \log{\frac{\mu_4}{\mu_3}}\,,\quad \centerG_2=\frac{C_2}{6 l} \log{\frac{\mu_2^2}{\mu_3\mu_4}}.$$
Finally, $C_1= \mu_1 e^{\frac{\centerG_1^2+\centerG_2^2}{C_2}}.$
We have shown that, provided $\{(\mu_1,\mu_2,\mu_3,\mu_4)\in\realpositive^4\ :\ \mu_2\mu_3\mu_4<\mu_1^3\}$, there is only one compatible Gaussian function, whose parameters are stated in the thesis.
\end{proof}

\subsection{Estimation errors}\label{sec:errors}
Let us assume that the measurements be affected by additive noise. This section is devoted to study the error committed computing the parameters from such noisy measurements.

For $i\in\until{4}$, let $\Delta_i$ be the error affecting measurement $\mu_i$ and let $\Delta$ be the column vector having components $\Delta_i$. If we recall the definition \eqref{eq:def-Phi} of the function $\Phi$, then the first order approximation of the error committed about the estimated parameters of $F$ can be computed as $(D\Phi)^{-1} \Delta$. Let us assume that the measurement errors $\Delta_i$ be independent Gaussian random variables with $\Exp[\Delta_i]=0$ and common variance $\Exp[\Delta_i^2]=\sigma^2$. Then the first order approximations of the errors on the estimates of the parameters of $F$ are Gaussian random variables: they have zero mean and their variance can be computed, after some algebra, as
\begin{align*}
S(l; C_2)=&\sigma^2 \frac{C_2^4 e^{2\frac{|\centerG|}{C_2}}}{9 C_1^2 l^4}\Big(9+e^{2\frac{l (l-2 \centerG_2)}{C_2}}+e^{2\frac{l \left(l+\sqrt{3} \centerG_1+\centerG_2\right)}{C_2}}
%\\&
+e^{2\frac{l \left(l-\sqrt{3} \centerG_1+\centerG_2\right)}{C_2}}\Big);
\end{align*}

\begin{align*}
S(l;C_1)=&\sigma^2 \frac{9 e^{2\frac{|\centerG|^2}{C_2}}}{l^4}
\Big(9 (|\centerG|^2-l^2)^2+(|\centerG|^2+2l\centerG_2)^2e^{\frac{2 l (l-2 \centerG_2)}{C_2}}\nn\\
&+\big(|\centerG|^2-l(\sqrt{3} \centerG_1+\centerG_2)\big)^2e^{\frac{2 l \left(l+\sqrt{3} \centerG_1 + \centerG_2\right)}{C_2}}\nn\\
&+\big(|\centerG|^2+l(\sqrt{3} \centerG_1-\centerG_2)\big)^2e^{\frac{2 l \left(l-\sqrt{3} \centerG_1+\centerG_2\right)}{C_2}}\Big)\,;
\end{align*}

\begin{align*}
S(l;|\centerG|)=&\sigma^2 \frac{C_2^2 e^{2\frac{|\centerG|^2}{C_2}}}{36 C_1^2 l^4 |\centerG|^2}
\Big(36|\centerG|^4+(2|\centerG|^2+2l\centerG_2)^2e^{\frac{2 l (l-2 \centerG_2)}{C_2}}\nn\\
&+\big(2|\centerG|^2+l(-\sqrt{3} \centerG_1-\centerG_2)\big)^2e^{\frac{2 l \left(l+\sqrt{3} \centerG_1 + \centerG_2\right)}{C_2}}\nn\\
&+\big(2|\centerG|^2+l(\sqrt{3} \centerG_1-\centerG_2)\big)^2e^{\frac{2 l \left(l-\sqrt{3} \centerG_1+\centerG_2\right)}{C_2}}\Big)\,,\nn
\end{align*}
where $|\centerG|=\sqrt{\centerG_1^2+\centerG_2^2}$;

\begin{align*}
S(l;\arctan{\frac{\centerG_2}{\centerG_1}})=&\sigma^2 \frac{C_2^2 e^{2\frac{|\centerG|^2}{C_2}}}{36 C_1^2 l^2 |\centerG|^4}\Big(4\centerG_1^2 e^{\frac{2 l (l-2 \centerG_2)}{C_2}}\\
&+(\centerG_1-\sqrt{3} \centerG_2)^2 e^{\frac{2 l \left(l+\sqrt{3} \centerG_1 +\centerG_2\right)}{C_2}} \nn \\
&+(\centerG_1+\sqrt{3} \centerG_2)^2 e^{\frac{2 l \left(l-\sqrt{3} \centerG_1+\centerG_2\right)}{C_2}}\nn
\Big).
\end{align*}

\begin{remark}[Center position and errors]\label{rem:center&errors}
Remark that all the above functions grow unbounded as $|\centerG|$ goes to infinity. This is consistent with the intuition that the sensors be more effective if close to the center of the Gaussian. Moreover, $S(l;\arctan{\frac{\centerG_2}{\centerG_1}})$ and $S(l;|\centerG|)$ grow unbounded also when $|\centerG|\to 0$. This suggests the opportunity of having nodes enough far apart from each other.
\end{remark}

Hence, there is some non-trivial optimization issue about $l$. The problem of a local optimization of the inter-node distance is the topic of the next paragraph.

\subsection{Optimal spacing}
Since we have restricted ourselves to hexagonal regular networks, the design parameter we are left with for optimization is $l$. It is clear that the choice of $l$ is constrained by the area to be covered and the number of available sensors.

Although we are aware of these potential constraints, in this section we keep them aside, and we we consider the problem of finding, for a fixed node, the value of $l$ which minimizes each of the above errors variances. That is, we address the unconstrained optimization problem
\begin{equation}\label{eq:lopt}
\lopt=\argmin_{l>0}{S(l)}.
\end{equation}
Whenever it is useful to specify which function we are referring to, we shall write $\lopt(C_1)$, $\lopt(C_2)$, and so on.

\begin{remark}[Finding $\lopt$]\label{rem:lopt}
Some facts on the optimization problem \eqref{eq:lopt} are immediate.
\begin{enumerate}
\item The errors variances $S(l;\cdot)$ depend continuously on $l$ and grow unbounded as $l$ goes to $0$ or to $\infty$. Hence there is at least one optimal $l$. Finding closed form expressions for these minima seems an intractable problem, but the minima can be evaluated numerically.
\item The parameter $C_1$ has no effect on any $\lopt$.
\item The parameter $C_2$ can be set as $C_2=1$ without loss of generality, since the mapping $(\centerG, C_2)\mapsto (\lambda \centerG, \lambda^2 C_2)$ induces a mapping as $\lopt \to \lambda \lopt$, for all $\lambda\in\realpositive.$
\end{enumerate}
\end{remark}
Hence $\centerG$, the relative location of the center, is the most meaningful variable. Before showing some numerical results, we state a bound about $S(l;C_2)$: %its tightness is numerically demonstrated in Figure~\ref{fig:lopt-stima}, while
its proof, which is based on direct computations, is postponed to the Appendix.
\begin{proposition}\label{prop:boundsC2}
Consider the optimization problem \eqref{eq:lopt}. Then
\begin{equation*}
\sqrt{C_2}-|\centerG|<\lopt(C_2)<\sqrt{2}\sqrt{C_2}+|\centerG|.
\end{equation*}
Namely, if $\centerG=0$, then $\lopt\cong1.11691\sqrt{C_2}.$
\end{proposition}

\begin{figure}
\begin{center}
\includegraphics[width=.7\columnwidth]{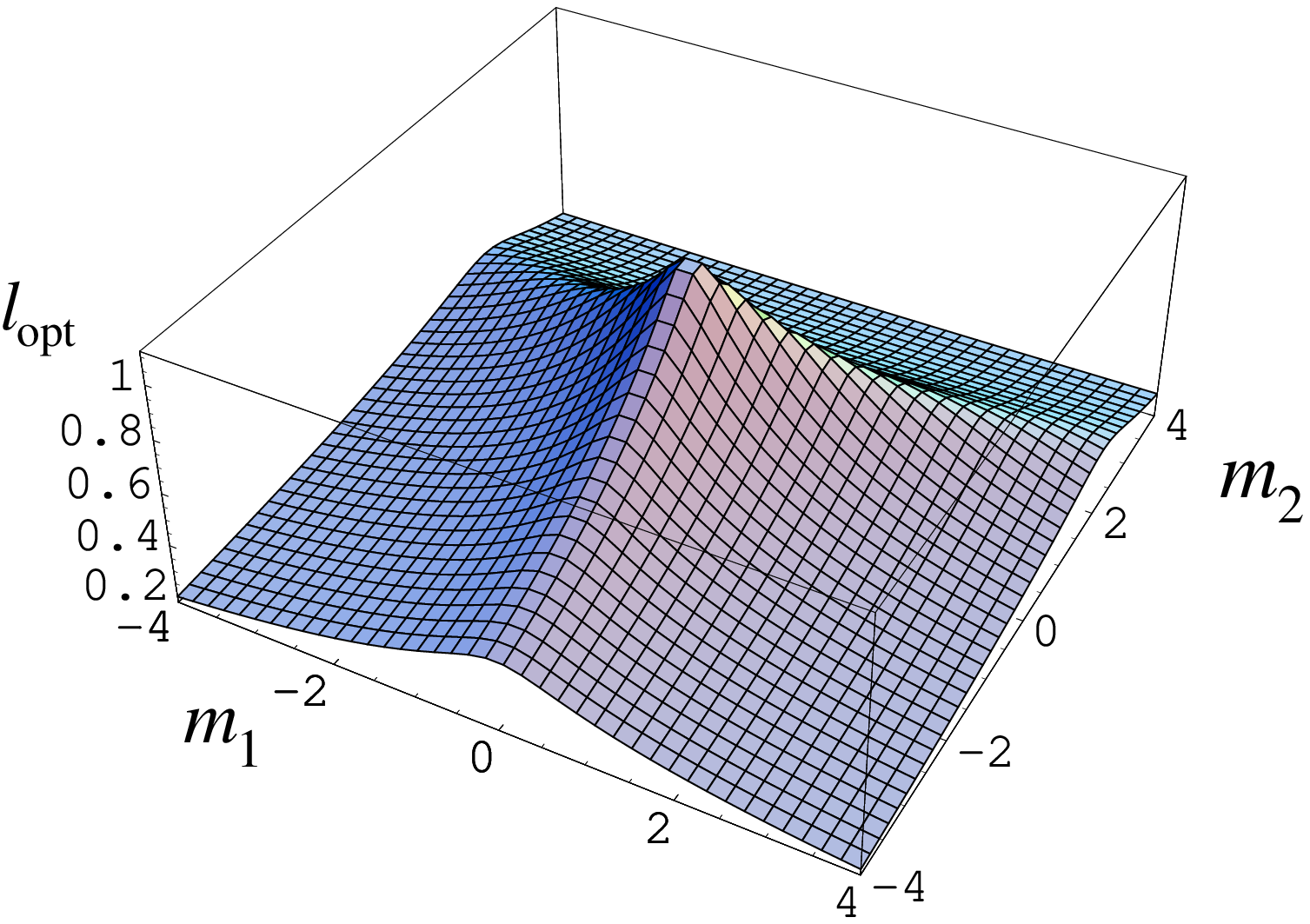}
\caption{Numerical computation of $\lopt(C_2)$, depending on $\centerG_1, \centerG_2$, with ${C_2=1}$.}
\label{fig:lopt2}
\end{center}
\end{figure}

\begin{figure}
\begin{center}
\includegraphics[width=.7\columnwidth]{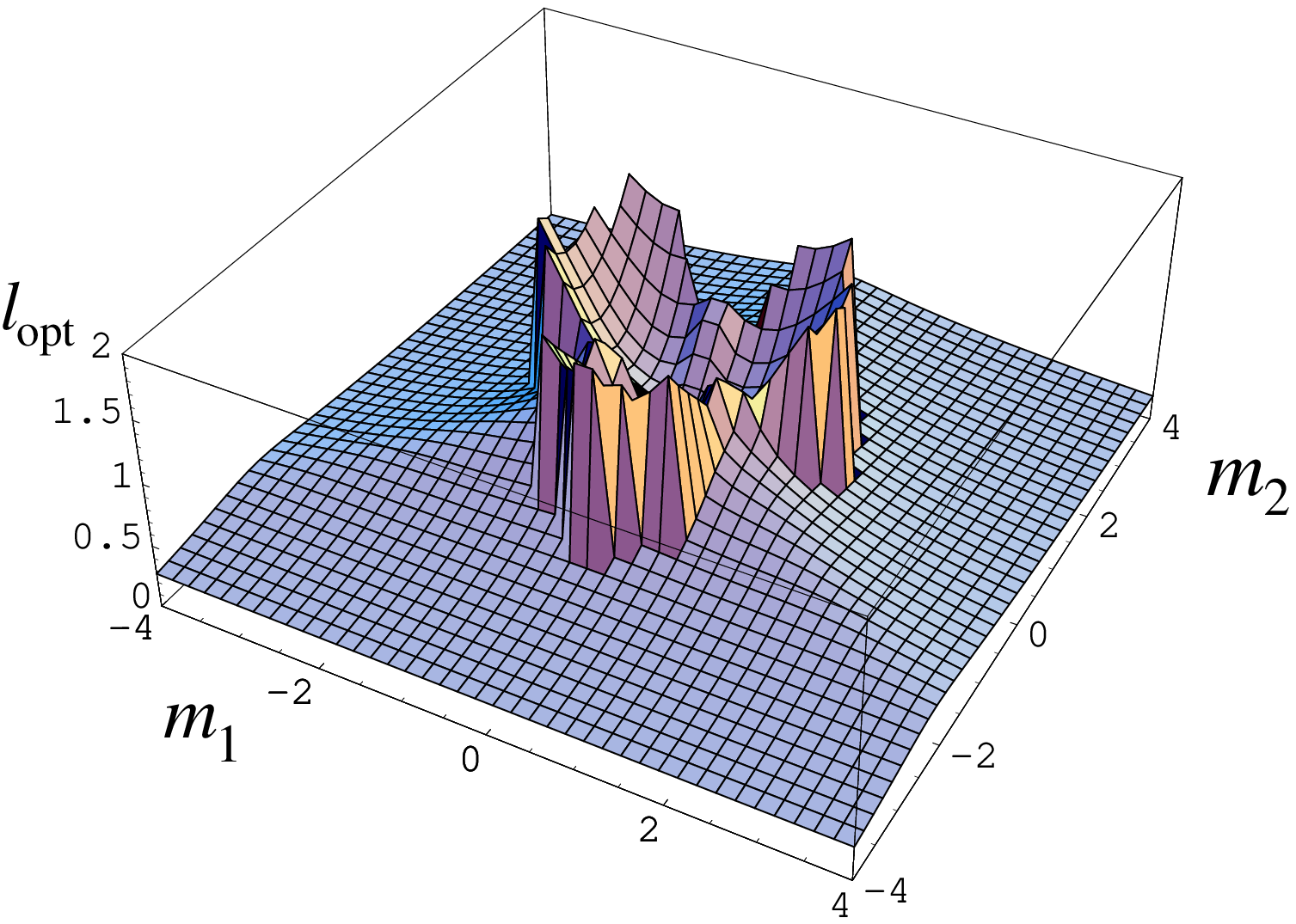}
\caption{Numerical computation of $\lopt(|\centerG|)$, depending on $\centerG_1, \centerG_2$, with ${C_2=1}$.}
\label{fig:lopt3}
\end{center}
\end{figure}

In Figures~\ref{fig:lopt2} and~\ref{fig:lopt3}, %~\ref{fig:lopt-cfr},~\ref{fig:lopt4},~and~\ref{fig:lopt5},
we show some computational results about the minima. %, obtained using the software {\textsf Mathematica.}
Some remarks are in order.
\begin{enumerate}
\item The optimal $\lopt$ depends on the parameters value, and namely on the value of $\centerG$ with respect to the central node, which is assumed to be located at the origin. As a consequence, the performed optimization is node-dependent, and can be of interest for off-line design when it is available some {\it a priori} information about the parameters.
\item The functions $S(l)$ have several local minima, and this results in discontinuities of the computed global minimum point $\lopt$.
\item It appears that $\lim_{|\centerG|\to \infty}\lopt=0.$ However, we know from Remark~\ref{rem:center&errors} that in such limit $S(l)$ diverges, for any $l$. Hence the estimate becomes unreliable, and the optimization looses its significance. We shall see, in the next section, how a data fusion algorithm can effectively cope with this problem, focusing on sensors close to the Gaussian center, so that the optimization of $l$ for far-away sensors becomes irrelevant.
\end{enumerate}

%\begin{figure}
%\begin{center}
%\includegraphics[width=.7\columnwidth]{stima-basso-mod.eps}
%\caption{Comparing the actual value of $\lopt(C_2)$ and the lower bound \eqref{eq:LBfigure}, assuming $\centerG_2=-1$ and $C_2=1$.}
%\label{fig:lopt-stima}
%\end{center}
%\end{figure}

%%%%%%%%%%%%%%%%%%%%%%%%%%%%%%%%%%%%%%%%%%%%%%%%%%%%%%%%%%%%%%%%%%%%%%%%%%%%%%%%%%%%%%%%%%%%%%%%%%%%%%%%%%%%%%%
\section{Fusion by consensus}\label{sec:Fusion}
In this section we discuss how the estimates of the parameters, obtained locally by the single nodes, can be fused at the network level into an agreement about the parameters of the Gaussian field. To the goal of fusion, we can take advantage of the wide literature \cite{ROS-JAF-RMM:07} on distributed consensus problems. Such problems can be solved by iterative distributed algorithms as the one we describe.

To work with more generality, let us consider a quantity of interest $x\in \reals$, which in the case of this paper will be one of the parameters of the Gaussian $F$. The algorithm evolves through discrete time steps $t \in \nonnegativeinteger$. The scalar $x_i(t)$ denotes the local estimate that, at time $t$, sensor $i\in \until{N}$ has of the quantity $x$: the estimate obtained from the estimation phase, described in Section~\ref{sec:Estimation}, will be $x_i(0)$. Note that, in our application, only the inner nodes of the network do possess initial estimates, and then take part to the fusion process. Thus $N$ is the number of inner nodes in the network. The iteration of the algorithm consists in
\begin{align}\label{eq:ConsensusAlgo}
x(t+1)=P(t) x(t),
\end{align}
where $x(t)=[x_1(t),\ldots, x_N(t)]^T$ is the estimates vector, and $P(t)$ is a stochastic matrix, depending on the {\em network topology} and on time. Namely, $P_{ij}(t)=0$ if $i$ and $j$ are not neighbors. The algorithm above, under suitable assumptions \cite{RC-FF-AS-SZ:08}, can be proved to converge to a consensus, in the sense that for all $i\in\until{N}$, we have $\lim_{t\to \infty}x_i(t)=\bar x$. The value $\bar x$ is said to be the consensus value, and it is a convex combination of the initial values $x(0)$. In some cases, the weights of this convex combination can be chosen by the designer. For instance, if $P(t)$ is constantly equal to a doubly stochastic matrix $P$, then $\bar x$ is the arithmetic average $\xave=N^{-1}\sum_{i=1}^{N}x_i(0)$.

Averaging is the first natural idea for fusion of the estimates obtained by each node. Indeed, elementary probability formulas imply that, given that each estimate $x_i(0)$ be affected by an independent random additive error whose variance is $\sigma^2_i$, the average $\xave=N^{-1}\sum_{i=1}^{N}x_i(0)$ is affected by an error of variance $N^{-2}\sum_{i=1}^{N}\sigma^2_i.$ In facts, if all the errors are equal, this yields reducing the error by a factor $N$. However, this is not the case of the present paper, when the {\em quality} of the estimates varies among nodes. Indeed, it can very well happen that, for some node $j$, the error variance be $\sigma^2_j < N^{-2}\sum_{i=1}^{N}\sigma^2_i.$ We argue that computing the arithmetic average is not always the best choice.
To cope with this problem we need to use a better estimator than the average.
Given the noisy estimates $x_i(0)$, such that $\Exp[(x_i(0)-x)^2]=\sigma_i^2$, the optimal estimator of $x$ in a mean squared error sense is
\begin{equation*}
\hat x=\frac{1}{\sum_{k=1}^{N} \frac{1}{\sigma_k^2} } \sum_{i=1}^{N}\frac{1}{\sigma_i^2} x_i(0).
\end{equation*}
If the variances $\sigma^2_k$ for all $k\in \until{N}$ were known {\it a priori}, then the problem could be solved as follows, combining two consensus algorithms. For every node $k\in\until{N}$, let $a_k(0)=\frac{x_i(0)}{\sigma_i^2},$ and $b_k(0)=\frac{1}{\sigma_i^2}.$ Then, the time-dependent vectors $a(t)$ and $b(t)$ evolve following two standard average consensus algorithms. We obtain that, for all $j\in\until{N}$,
$$\lim_{t\to\infty}\frac{a_j(t)}{b_j(t)}=
\frac{\frac{1}{N}\sum_{l=1}^{N} \frac{x_l(0)}{\sigma^2_l}}{\frac{1}{N}\sum_{m=1}^{N} \frac{1}{\sigma^2_m}}
=\hat x.$$

Instead, the variances $\{\sigma_k\}$ are not {\it a priori} known in the case of this paper, because they depend on the value $x$, as illustrated in Section~\ref{sec:errors}. For this reason, the above procedure can not be used, and we propose instead a suboptimal distributed algorithm.
At each time step, every node $i\in\until{N}$ updates and communicates both the variable $x_i(t)$, and a variable $s_i(t)\in \realpositive,$ which approximates the variance of the error that $i$ is committing in estimating $x$ with the value $x_i(t)$. The algorithm we propose is as follows.

Let us consider the graph $\G$ of the connections among the inner nodes of the hexagonal network, whose nodes are labeled in $\until{N}$. For all $i\in \until{N}$, let $\neigh_i\subset \until{N} $ denote the set of the neighbors of $i$ in the graph, including $i$ itself. Let $x(0)$ be the vector of the local estimates, and $s(0)$ the vector of the presumed variances of the errors affecting $x(0)$, as computed by each node $i$ substituting $x_i(0)$ into the formulas in Section~\ref{sec:errors}.

At each time step $t\in\nonnegativeinteger$, each node $j$ communicates its states $x_j(t)$ and $s_j(t)$ to all its neighbors, and the states are updated following
\begin{align}\label{eq:wise}
x(t+1)&=P(t) x(t)\\
\nonumber s(t+1)&=M(t) s(t),
\end{align}
where
$$P_{ij}(t)=\begin{cases}\displaystyle\frac{1/s_j(t)}{\sum_{k\in \neigh_i}{1/s_k(t)}}, & \text{ if $j\in \neigh_i$,} \\
\qquad 0 & \text{ otherwise};
\end{cases}$$
and
$$M_{ij}(t)=\begin{cases}\displaystyle \frac{1/(s_j(t))^2}{\sum_{k\in \neigh_i} 1/(s_k(t))^2}, & \text{ if $j\in \neigh_i$,} \\
\qquad 0 & \text{ otherwise}.
\end{cases}$$
%\begin{align*}
%M_{ij}(t)&=\frac{(P_{ij}(t))^2}{\sum_{k\in \neigh_i}{(P_{ik}(t))^2}}\\
%&=\frac{1/(s_j(t))^2}{\sum_{k\in \neigh_i} 1/(s_k(t))^2},
%\end{align*}
%if $j\in \neigh_i$, and zero otherwise.

The above update can be done in a distributed way, using information at the node and received from the immediate neighbors. Remark that the update matrices depend on the state of the system, and then on time. The convergence of the proposed algorithm is stated in the following result.

\begin{proposition}[Convergence]
Provided the graph $\G$ is connected, system~\eqref{eq:wise} converges to consensus in both $x(t)$ and $s(t)$, that is, there exist scalars $x^*,s^*$ such that
$$\lim_{t\to \infty}
\left(
  \begin{array}{c}
  x(t) \\
  s(t) \\
  \end{array}
\right)=
\left(
  \begin{array}{c}
    x^*\1 \\
    s^*\1 \\
  \end{array}
\right),$$
where $\1$ is the $N\times 1$ vector whose entries equal 1.
\end{proposition}
\begin{proof}
Unlike time-independent consensus iterations \cite{RC-FF-AS-SZ:08}, an eigenvalue analysis is not in general sufficient to conclude convergence for a time-varying system. Instead, we take advantage of the very general convergence results in \cite{LM:05}. To this goal, we have to check that the linear iterations~\eqref{eq:wise} satisfy the Strict Convexity Assumption in \cite[page~172]{LM:05}. We make the following remark: every $M(t)$ is a stochastic matrix, and then $s_j(t)\in [\min_i s_i(0), \max_i s_i(0)]$, for all $t\ge 0$, and all $j\in \until{N}$. Since these intervals are compact and do not contain zero, we argue that for all $t\ge0$, the non-zero entries of $P(t)$ and $M(t)$ belong as well to a compact set of positive numbers. This is enough to satisfy the convexity assumption, and since $\G$ is connected, we can apply \cite[Th.~3]{LM:05} and conclude the proof.
\end{proof}

Note that, while the state space of each agent is two-dimensional, the object of interest is the consensus about the variable $x(t)$ only. The consensus value $x^*$ is not in general equal to the optimal estimator $\hat x$, but the simulations presented in Section~\ref{sec:Simul} show that $x^*$ is much closer to the true value $x$, than the arithmetic average $\xave$.

\begin{remark}[Algorithm variations]
The proposed algorithm is not the only heuristic solution: this remark informally illustrates two other possibilities. A very simple idea is the following one. Instead of running a consensus algorithm about $s(t)$, the variables $s_i(t)$ can be computed at each time step from the current estimates $x_i(t)$, using the formulas in Section~\ref{sec:errors}. A combination of this latter approach with a consensus scheme can also be pursued, as follows. At each time step $t$, each node $i$ computes from $x_i(t)$ an estimate $s_i(t;k=0)$, and then runs some (say $\bar k$) consensus steps about $s$ only. The obtained $s_i(t;\bar k)$ is then used in a consensus step to compute $x_i(t+1)$. A convergence proof for these variations can be adapted from the one given above.
\end{remark}

\begin{figure}[htb]
\begin{center}
\includegraphics[width=0.49\textwidth]{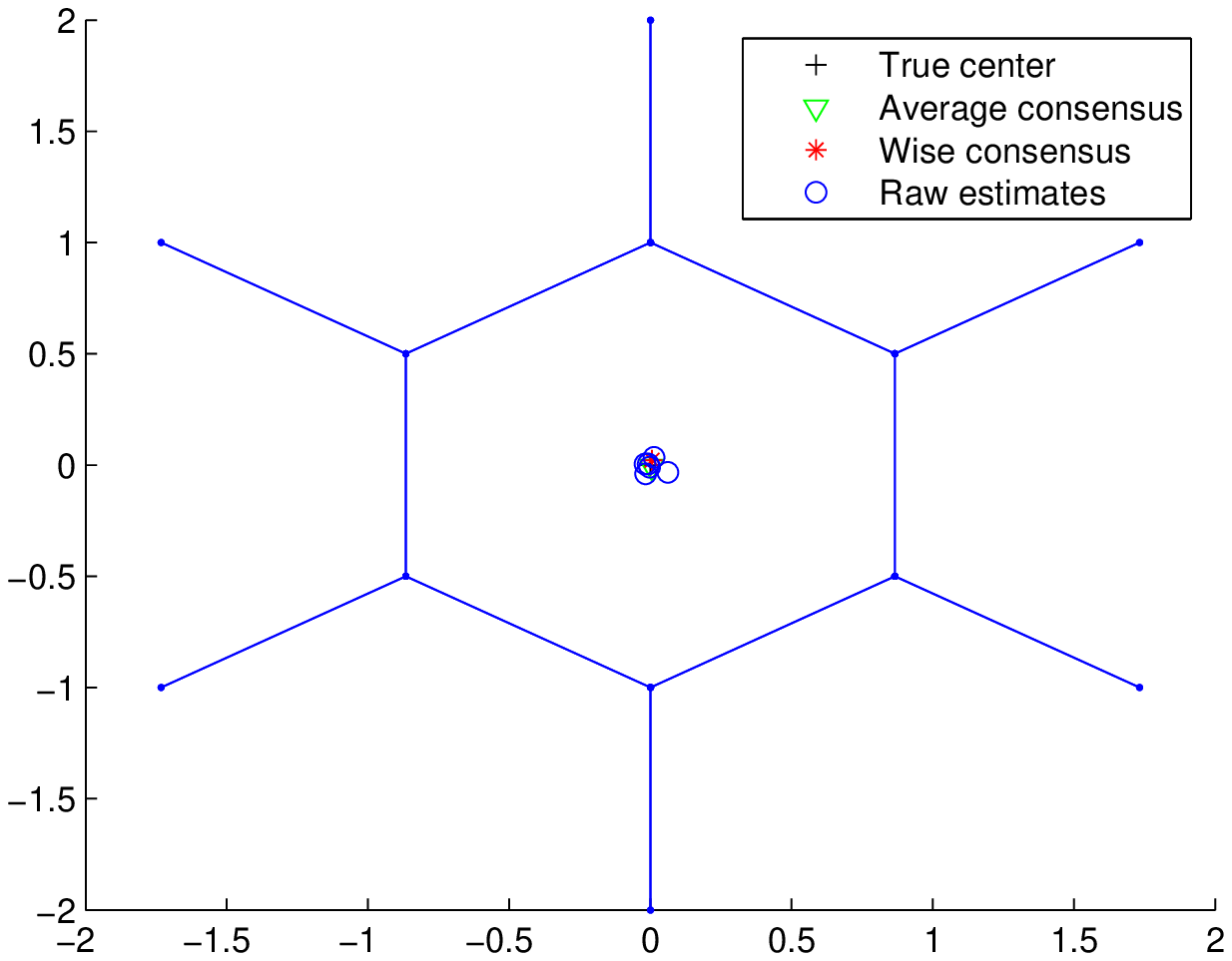}
\includegraphics[width=0.49\textwidth]{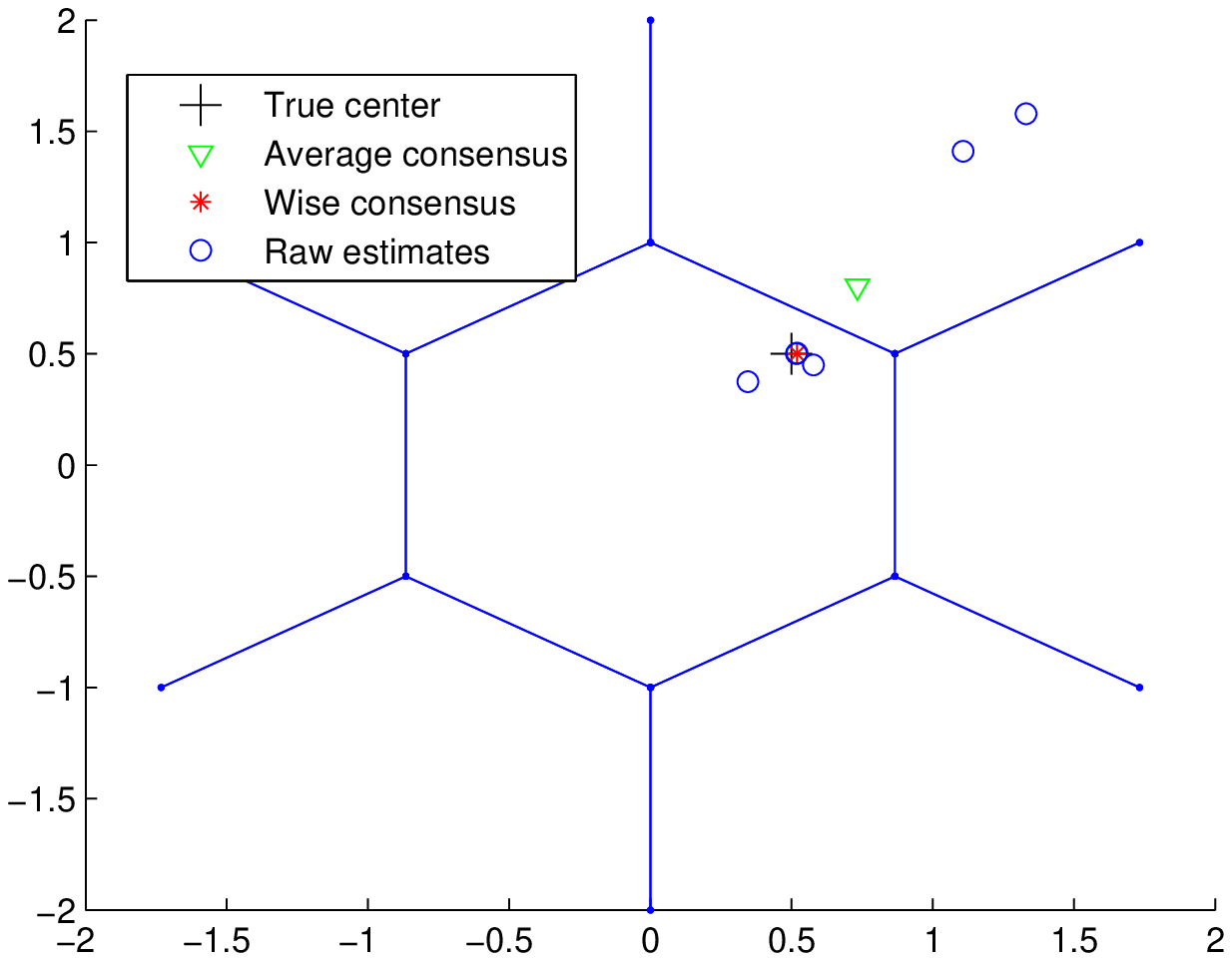}
\includegraphics[width=0.49\textwidth]{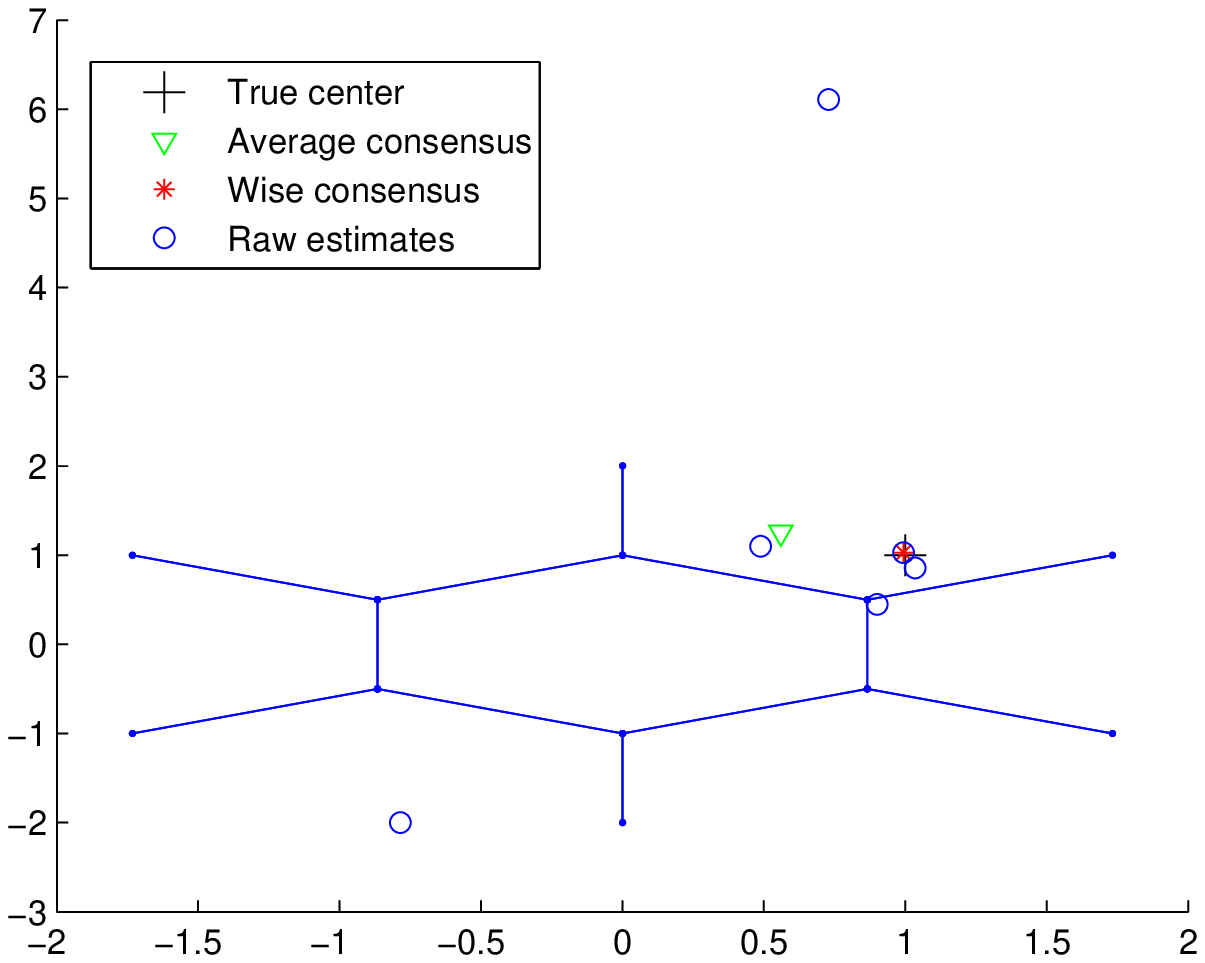}
\includegraphics[width=0.49\textwidth]{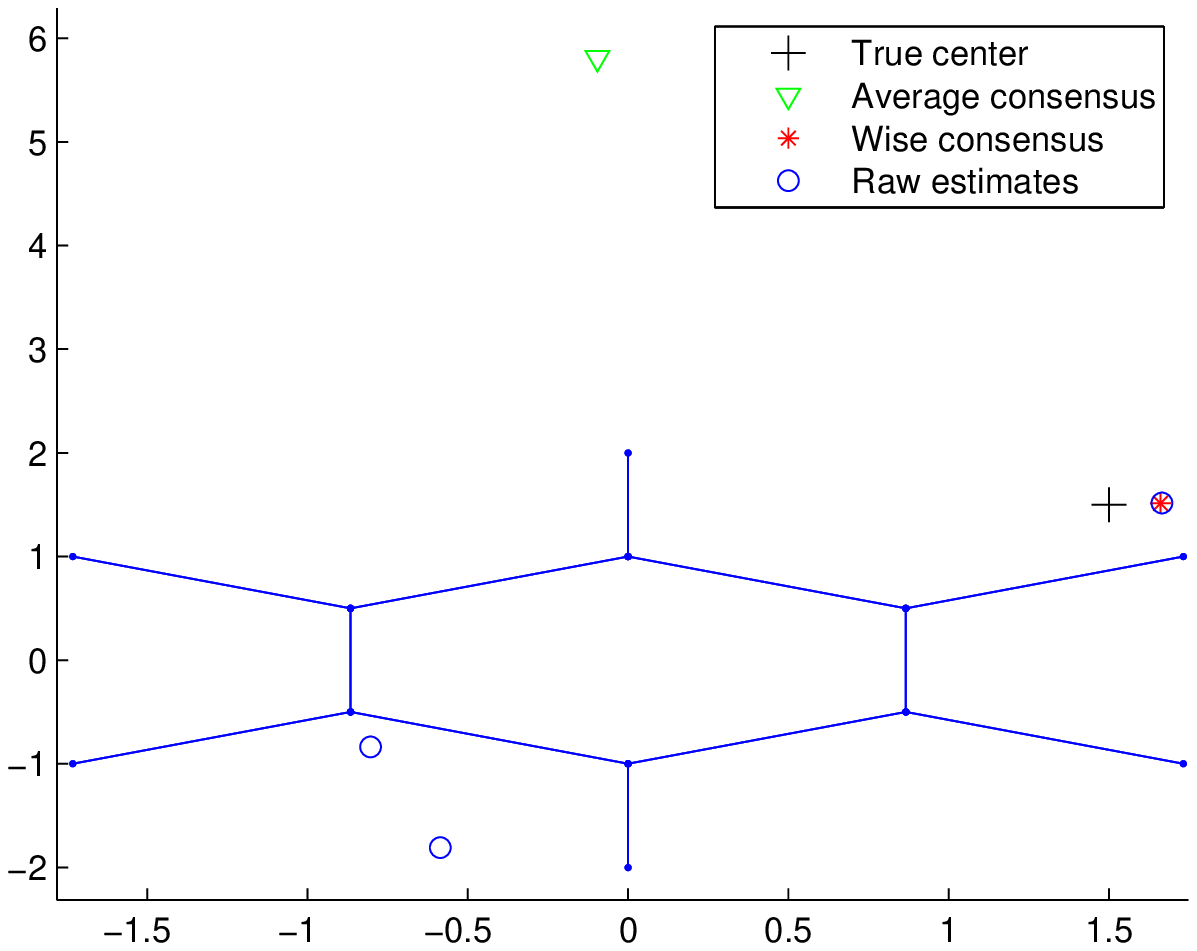}
\caption{Raw and fused estimates of the center of the Gaussian $F$, whose position $(\centerG_1,\centerG_2)$ is equal to $(0,0)$, $(0.5,0.5)$, $(1,1)$, and $(1.5,1.5),$ respectively. The label {\it Wise consensus} denotes the algorithm \eqref{eq:wise}.}
\label{fig:HexEst}
\end{center}
\end{figure}
\section{Simulations}\label{sec:Simul}
Simulations are presented to demonstrate the effectiveness of the proposed estimation and fusion algorithm. To provide an example, we focus on the estimation of the center $(\centerG_1, \centerG_2)$ of the Gaussian field. We consider a network of twelve nodes, six of which are inner nodes forming an hexagon, and hence are able to compute estimates and participate to the fusion process. We assume that each measurement is corrupted by an additive random Gaussian variable whose variance is $1\%$ of the maximum value of the signal to be estimated. Figure~\ref{fig:HexEst} shows the raw estimates obtained by each node, and the fused estimate of the center, compared with its true position. Note that, although the measurement error seems to be low, it is enough to perturb significantly the local estimates. Hence, the simulations setting is suitable to demonstrate the effectiveness of the proposed fusion algorithm.

Simulations show that the estimation of the bi-dimensional parameter is possible with the proposed distributed method, and that the sensitivity to measurement errors is greatly reduced by fusion: the consensus algorithm proposed in \eqref{eq:wise} is much more accurate, and then preferable to simple averaging. As expected, the sensor network is less effective when the Gaussian is eccentric with respect to the network. Nevertheless, for moderate eccentricity the fusion algorithm gives an effective compensation of the errors affecting the local estimates.

\section{Conclusion and open problems}\label{sec:outro}
In this work we have posed and solved the problem of distributely estimating a Gaussian signal from noisy measurements taken by networked sensors. Our solution consists in a two-phases algorithm designed for a network which forms an hexagonal tessellation. In a first phase, sensors locally compute estimates of the parameters identifying the Gaussian signal, using their own measurement, and the measurements by their immediate neighbors. In the second phase, an heuristic consensus algorithm is run to obtain a common estimate. The consensus phase takes advantage of an estimation of the statistics of the errors affecting the local estimates. This feature allows to obtain encouraging results in terms of accuracy, also in presence of significant measurement noise, which yields large errors in the local estimates. Besides proving the convergence of the algorithm, we also studied the reliability of the local estimates, obtaining useful results for the network design. Our work is suitable of several further developments, that we briefly outline as follows.

\begin{enumerate}
\item The analysis given in the present paper can be extended, through similar computations, to the case of a more general function $$F(x)=C_1 e^{-\frac{(x-m)^T\Gamma (x-m)}{C_2}},$$ provided $\Gamma$ is a known positive definite matrix. This version of the problem has a clear interpretation in terms of the application to wild-fires monitoring, whenever there is an anisotropy in the fire spreading (e.g. due to prevalent winds). Beyond this extension, an approach similar to the one of the present paper is likely to be useful on more general parametric estimation problems.

\item Assuming that the sensors be deployed to form an hexagonal network offers some advantages, but is also a restriction, compared to general sensor networks. This remark suggests some natural developments. First, a necessary step towards application is a sensitivity analysis with respect to errors in the deployment of sensors. Second, the application of the algorithm would benefit from the complementary design of distributed algorithms to achieve the hexagonal network by deployment of self-propelled sensors. Third, there can be cases in which the hexagonal configuration is difficult or impossible to achieve. Hence, it would be useful to extend our analysis to different, possibly non regular, graphs. Natural questions include: which graphs are suitable to solve the Gaussian parameters estimation problem? Among those, which ones are preferable? After these basic questions are answered, a locational optimization problem could be defined and hopefully solved in a distributed fashion, to achieve an optimal sensors deployment. Note that problems of distributed locational optimization have already been approached in the literature, for instance in \cite{JC-SM-TK-FB:02j,FB-JC-SM:09}.

\item The literature on distributed sensor fusion has already devoted a considerable attention to the design of algorithms to compute the average of local values. As we have shown, this can happen to be a non-optimal solution to a fusion problem, when local values have different (possibly unknown) reliability. Beyond any technical improvement of the proposed algorithm, we believe that further research should focus on the problem of an optimal distributed fusion of estimates having different ``quality''.
\end{enumerate}

\bibliographystyle{plain}
\bibliography{aliasFrasca,RefFrasca}

\appendix
\section{Proof of Proposition~\ref{prop:boundsC2}}
%\begin{proof}[of Proposition~\ref{prop:boundsC2}]
Up to a multiplicative constant, we can compute
\beqn
\frac{d}{dl} S (l)&=&\frac{4}{l^5}\Big(-9+\frac{l^2-l \centerG_2-C_2}{C_2}e^{2\frac{l (l-2 \centerG_2)}{C_2}}\nn\\
&&\quad+\frac{l^2+\frac{\sqrt{3}}2 l\centerG_1+l\frac{\centerG_2}2 -C_2}{C_2}e^{2\frac{l \left(l+\sqrt{3} \centerG_1+\centerG_2\right)}{C_2}}\nn\\
&&\quad+\frac{l^2-\frac{\sqrt{3}}2 l\centerG_1+l \frac{\centerG_2}2 -C_2}{C_2}e^{2\frac{l \left(l-\sqrt{3} \centerG_1+\centerG_2\right)}{C_2}}\Big).\nn
\eeqn
Hence, if $\centerG=0,$ then the necessary condition reads $\left(\frac{l^2}{C_2}-1\right) e^{2\frac{l^2}{C_2}}=3,$ which can be numerically solved, giving $l\cong 1.11691 \sqrt{C_2}$.
If instead $m\neq0$, then the minimum has to satisfy at least one of the following inequalities
\begin{eqnarray}
l^2-l \centerG_2-C_2&>&0\nn\\
\quad l^2+\frac{\sqrt{3}}2 l\centerG_1+l\frac{\centerG_2}2 -C_2&>&0\label{disug}\\
\quad l^2-\frac{\sqrt{3}}2 l\centerG_1+l \frac{\centerG_2}2 -C_2&>&0\nn
\end{eqnarray}
That is,
\begin{align*}
\lopt>\min&\left\{\frac12\Big(\centerG_2+\sqrt{4 C_2+\centerG_2^2}\,\Big),\right.
\\&
\left.\frac12\Big(-\frac{\centerG_2}2 -\frac{\sqrt{3}\centerG_1}2 +\sqrt{4 C_2+\Big(\frac{\centerG_2}2 +\frac{\sqrt{3}\centerG_1}2\Big)^2}\,\Big),\right.\\
&\left.\frac12\Big(-\frac{\centerG_2}2 +\frac{\sqrt{3}\centerG_1}2+\sqrt{4 C_2+\Big(\frac{\centerG_2}2 -\frac{\sqrt{3}\centerG_1}2\Big)^2}\,\Big)\right\}.
\end{align*}
Such a minimum is always positive, and moreover
\begin{eqnarray}
\nonumber\lopt&>&\frac12\Big(-|\centerG|+\sqrt{4 C_2+|\centerG|^2}\Big)\\%\label{eq:LBfigure}\\
\nonumber&=&\frac12 \frac{4 C_2}{|\centerG|+\sqrt{4 C_2+|\centerG|^2}}\\
\nonumber&>&\frac{C_2}{|\centerG|+\sqrt{C_2}}\\
\nonumber&>&\sqrt{C_2}-|\centerG|,
\end{eqnarray}
using the Cauchy-Schwartz inequality and the fact that the function $y\mapsto y+\sqrt{4 C_2+y^2}$ is increasing.

Let us now look for an upper bound on $\lopt$. To this goal, let us assume $l$ to satisfy
\beqn
l^2-l |\centerG| -C_2>0\,, \label{coeff-pos}
\eeqn
so that the inequalities~\eqref{disug} are all satisfied. Then, using $e^y>1+y$, we argue that
\beqn
\frac{d}{dl} S (l)&>&l^4+l^2 (|\centerG|^2-C_2/2)-2 C_2^2. \label{stima-raff}
\eeqn
Actually the right-hand-side of \eqref{stima-raff} is positive when
$$l^2>\frac12 \Big(C_2/2-|\centerG|^2+\sqrt{(C_2/2-|\centerG|^2)^2+8 C_2^2}\,\Big).$$
Since
\begin{align*}
\sqrt{\frac12 \Big(C_2/2-|\centerG|^2+\sqrt{(C_2/2-|\centerG|^2)^2+8 C_2^2}\Big)}&<\sqrt{|C_2/2-|\centerG|^2|+\sqrt{2}C_2}\nn\\
&<|\centerG|+\sqrt{2C_2}
\end{align*}
and, for every $l> |\centerG|+\sqrt{2C_2}$, Equation~\eqref{coeff-pos} is satisfied, we obtain the upper bound
$$\lopt< |\centerG|+\sqrt{2C_2}\,.$$

%\end{proof}

\end{document}